\documentclass{article}

\usepackage{microtype}
\usepackage{graphicx}
\usepackage{subfigure}
\usepackage{xcolor}
\usepackage{amsmath}
\numberwithin{equation}{section}
\usepackage{booktabs} 

\usepackage{hyperref}

\usepackage[accepted]{icml2025}

\usepackage{amsmath}
\usepackage{amssymb}
\usepackage{mathtools}
\usepackage{amsthm}

\usepackage[capitalize,noabbrev]{cleveref}

\theoremstyle{plain}

\theoremstyle{definition}

\theoremstyle{remark}

\usepackage[textsize=tiny]{todonotes}

\icmltitlerunning{OPTIMUS: Optimization Productivity Tool for Intelligent Management of Utilizable Space}

\begin{document}

\twocolumn[
\icmltitle{OPTIMUS: Optimization Productivity Tool for Intelligent Management of Utilizable Space}



\icmlsetsymbol{equal}{*}

\begin{icmlauthorlist}
\icmlauthor{Souvik Bhattacharyya,}{equal,comp}
\icmlauthor{Nisha Singh,}{equal,comp}
\icmlauthor{Salman Haider,}{comp}
\icmlauthor{Balaji Nagarajan,}{comp}
\icmlauthor{Ved Prakash Dwivedi,}{comp}
\icmlauthor{Nithin Surendran,}{comp}
\icmlauthor{Karthik Nair}{comp}
\end{icmlauthorlist}

\icmlaffiliation{comp}{Lowes', 1000 Lowe's Boulevard Mooresville, NC 28117 }

\icmlcorrespondingauthor{Ved Prakash Dwivedi}{vedprakash.dwivedi@lowes.com}

\icmlkeywords{Space optimization, POG optimization, Bay optimization, Dynamic programming, Linear Binary Knapsack problem}

\vskip 0.3in
]



\printAffiliationsAndNotice{\icmlEqualContribution} 

\begin{abstract}
We study department-level retail space optimization, where limited bay capacity must be allocated among planograms (POGs) under business and operational constraints. The problem is formulated as a linear binary knapsack model, with potential SKUs treated as items characterized by space requirements and weighted value contributions from sales, margin, units, and assortment similarity. Dynamic Programming (DP) is employed to obtain exact and reproducible assortment decisions in $O(nc)$ time, avoiding the variance inherent in heuristic approaches. These decisions are integrated with a second-stage bay optimization model formulated as a mixed-integer program. Evaluated end-to-end across ten optimization runs spanning multiple departments and store clusters, the OPTIMUS framework achieves an average sales lift of 11.8\% and an average margin lift of 9.5\%. Overall, OPTIMUS provides a scalable, interpretable, and profit-driven solution for enterprise-scale retail space management.
\end{abstract}

\section{Introduction}
In the retail industry, physical space is directly tied to revenue generation, making the efficient utilization of available space a critical challenge. Retail floor areas are typically divided into bays \cite{curhan1972relationship}, which represent the fundamental sellable units where products are displayed. These bays are then aggregated into planograms or POGs -- visual merchandising blueprints that define product placement and category layouts across stores. Planograms, in turn, roll upward into broader selling departments that drive profitability. It is at the department level that the problem of space optimization becomes most prominent, as it directly affects both store-level performance metrics and strategic decision-making \cite{corstjens1981model, dreze1994shelf}

Retailers operate under increasing pressure to maximize profitability within constrained store footprints. As customer demand for product variety grows, the tension between limited physical space and the need for broad assortments intensifies. This makes space allocation not only an operational task but also a strategic decision that influences category management, supplier negotiations, and overall store performance \cite{hubner2017integrated}.
Traditional research in this area has proposed a variety of models for space management \cite{urban1998inventory}, ranging from demand-based space elasticity frameworks \cite{dreze1994shelf} to heuristic algorithms and mixed-integer programming approaches \cite{van2010ordering}
While these methods offer valuable insights, they often operate at the bay or shelf level, focusing on micro-optimizations of product arrangements. However, such models can overlook the aggregated effects at the department level, where space allocation decisions yield the greatest managerial impact. Moreover, many existing formulations are nonlinear and computationally intensive, posing practical challenges for implementation in large-scale retail environments.
From an operations research perspective, the problem aligns closely with combinatorial optimization. The knapsack problem and its variants have emerged as useful abstractions for modeling resource-constrained allocation problems, including those found in retail \cite{kellerer2004multidimensional}. By adapting this framework to the retail context, it is possible to formalize the trade-off between limited shelf space and the maximization of profit in a mathematically tractable way.

Despite substantial research efforts, there remain two critical gaps in the literature. First, most existing studies focus on localized shelf or bay-level optimization without adequately addressing department-level decision-making. Second, many models rely on nonlinear demand functions or complex simulation methods that limit their applicability in practice. This creates a disconnect between academic approaches and the needs of retailers, who require scalable, computationally efficient solutions that can be readily integrated into operational workflows.
At the selling department level, where profitability and product strategy converge, the challenge is to allocate limited bay space in a way that maximizes returns while respecting category-level and organizational constraints. Given the increasing complexity of modern retail assortments, solving this optimization problem in a tractable and practical manner is essential.

This paper addresses the gap by framing department-level retail space optimization as a linearized binary knapsack problem. In this formulation, each product category represents an item with associated profit potential, while the department’s available bay space serves as the knapsack capacity. The binary structure models inclusion or exclusion decisions, and the linearization ensures computational feasibility for large-scale retail applications.
The contributions of this study are threefold:
Model Development – We introduce a novel application of the binary knapsack problem to department-level space allocation, bridging the gap between bay-level optimization and higher-level managerial decision-making.
Computational Tractability – We propose a linearized formulation that is solvable with standard optimization tools, making it accessible and scalable for real-world retail environments.
Managerial Insights – We demonstrate how the model can inform practical trade-offs between product variety, planogram design, and profitability, offering actionable guidance for category management and strategic planning.
By combining theoretical rigor with practical applicability, this work contributes to both the operations research literature and the retail industry, providing a framework for addressing one of the sector’s most pressing optimization challenges.

\section{Methodology}
To operationalize the proposed framework, we adopt a two-step optimization approach designed to reflect both the strategic and operational layers of retail space management. The first step, planogram optimization, determines which products should be selected for stocking on shelves. This decision is driven by key performance indicators such as historical sales volume and contribution margins, ensuring that the chosen product set maximizes expected profitability while maintaining assortment relevance. By framing this step as a selection problem, we effectively filter the wide product universe down to the most economically viable candidates for display.

Building on this foundation, the second step, bay optimization, addresses the allocation of physical shelf space across selling departments. Here, the focus shifts from item-level selection to space distribution, determining how much of the available bay capacity should be dedicated to each department based on the previously chosen assortments. This step balances the competing demands of different departments, ensuring that space is distributed in proportion to profitability potential while respecting organizational constraints. Together, the two steps create a structured pathway: first identifying the what (products to stock), and then optimizing the how much (space to allocate), thereby aligning product assortment strategies with space utilization for maximum overall returns.

\subsection{Planogram optimization}

A POG is defined as a structured set of items that together represent a specific product category or display unit within a retail store. At the most granular level, we have individual items (for example, different brands or models of refrigerators). When a collection of such related items is grouped together, it forms a planogram—for instance, a Refrigerator POG or a Dishwasher POG. Multiple planograms are then aggregated into a selling department, such as the Appliances Department, which may include refrigerator and dishwasher POGs.

To determine how much space should be allocated to a given POG, it is essential to estimate the value generated under different space configurations. The POG optimization problem is therefore defined as finding the optimal subset of items to include within a POG for a specified amount of available space \cite{gallego2019assortment}. This optimization is performed for every possible space allocation and for each POG. The resulting outputs then serve as inputs to the subsequent bay optimization stage, which determines the optimal space distribution among POGs within a department.

The total display area available for product placement within a planogram is defined as its capacity, which represents the upper limit of linear shelf space that can be allocated. Each POG contains a set of potential items—the SKUs that could be included in the assortment. These items may be local (currently stocked in the store and supported by actual sales data) or non-local (not presently carried but with sales potential estimated from comparable stores within the same regional cluster). Each item occupies a fixed amount of shelf space and is characterized by an average selling price (ASP) and an average margin per unit (AMPU). Because some products may exhibit substitution or complementarity effects, the overall demand associated with an assortment depends not only on the individual performance of items but also on their combined composition. These item-level parameters form the basis for the subsequent optimization, enabling the model to balance spatial constraints with revenue and profitability objectives.

To balance different business objectives, the optimization framework incorporates weighted preferences across multiple metrics—sales, margin, units sold, and assortment similarity. These weights are represented by the vector 

\begin{equation}
  \lambda = (\lambda_\text{sales}, \lambda_\text{margin}, \lambda_\text{units}, \lambda_\text{similarity})
\end{equation}
\label{objweights}

allowing users to emphasize specific goals, such as prioritizing margin over sales. Each item’s contribution to the objective function is computed as a weighted sum of these metrics.

We denote the subset of items included in a given assortment by an indicator vector 
$x \in \{ 0, 1 \}^n$, where $n$ is the total number of potential items. For an item $i$, the component $x_i=1$ indicates that the item is included in the assortment, while
$x_i=0$ indicates exclusion.

Using the inputs described above, we define the objective function

\begin{equation}
f \colon \{ 0, 1 \}^n \to [0, \infty)
\end{equation}
\label{objfunc}

and formalize the POG optimization problem as a variant of the binary knapsack problem. Two formulations are possible: a linear knapsack problem, which assumes no inter-item dependencies (i.e., no substitution or complementarity effects), and a quadratic knapsack problem, which accounts for pairwise item interactions. To simplify computation and enhance scalability, the current implementation employs the linear knapsack formulation, excluding direct item-to-item interaction effects. The consideration of quadratic effects is reserved for future extensions of this work.

\subsubsection{The Linear Binary Knapsack Problem}

In certain POGs, the interaction effects between items are negligible, unavailable, or statistically insignificant. When this occurs, it is reasonable to assume independence among the potential items within a POG. This assumption simplifies the overall optimization problem by removing substitution or complementarity effects, allowing it to be formulated as a linear binary knapsack problem. In this linear case, each item’s contribution to the total objective is constant and unaffected by the inclusion or exclusion of other items.

The linear objective function can be expressed as:
\begin{equation}
	f(x) = \sum_{i = 1}^n q_i x_i = x^\text{T} q
\end{equation}
where $n$ denotes the number of potential items, 
$x \in \{ 0, 1 \}^n$ represents the binary assortment vector, and $q$ is the \emph{objective vector}, with each $q_i$ indicating the contribution of item 
$i$ to the overall objective value.

The objective vector $q$ is derived component-wise using the weight vector as given in eq. \ref{objweights}, the estimated demand vector $d$, and the baseline assortment $x^0$, as follows

\begin{align}
	m_i & = p_i \lambda_\text{sales} + g_i \lambda_\text{margin} + \lambda_\text{units} \\
	q_i & = m_i d_i + (2 x^0_i - 1) \lambda_\text{similarity}
\end{align}
for $i = 1,\dots,n$. Here, the vector $m$ functions as the demand objective multiplier, translating expected demand into performance-based measures such as total sales, margin, or units sold. The similarity component, which is independent of demand, captures the deviation between the optimized and baseline assortments.
To compute the similarity term, $q_\text{similarity} = 2 x^0_i - 1$, the baseline assortment vector $x^0$ is transformed from 0 to a -1 using the mapping $x \to 2x - 1$. This ensures that both additions and removals of items ($0\to1$ and $1\to0$ transitions) are penalized symmetrically, thereby treating all potential items uniformly. Conceptually, this component can be interpreted as a constant (representing the total number of items in the baseline assortment) minus the Hamming distance between $x$ and $x^0$. The penalty for deviation is scaled by the weight 
$\lambda_{similarity}$ , providing flexibility in how strongly assortment changes are discouraged. When item-level operational data are available, this framework can be further refined by replacing the generic penalty with item-specific switching costs, offering a more accurate reflection of real-world merchandising constraints.

\subsubsection{POG optimization Algorithms}

A variety of algorithms can be employed to solve knapsack-style optimization problems \cite{ahmed2008solving}, each offering distinct trade-offs between computational efficiency, solution quality, and implementation complexity. Common approaches include Brute Force, Greedy, Greedy Productivity, Mixed Integer Programming (MIP), and Genetic Algorithms (GA). Brute Force methods guarantee optimality by evaluating all possible combinations but are computationally infeasible for large problem sizes. Greedy-based techniques are computationally efficient but often yield approximate solutions that may deviate from the global optimum. MIP formulations provide flexibility for incorporating complex constraints but can become computationally expensive as the problem scales. Metaheuristic methods such as GA introduce stochasticity, which enhances solution diversity but also leads to variability in outcomes and higher tuning costs.

In this study, we adopt Dynamic Programming as the primary optimization approach for solving the linear binary knapsack formulation of the POG optimization problem. It leverages the inherent optimal substructure and overlapping subproblems properties of the knapsack framework to compute exact solutions efficiently \cite{martello1990knapsack}. By decomposing the problem into smaller subproblems and systematically combining their solutions \cite{pisinger1999exact}, it achieves optimal results in pseudo-polynomial time—offering a balance between computational tractability and accuracy. Compared to brute force enumeration or heuristic techniques, dynamic programming provides a deterministic, reproducible, and computationally efficient foundation for POG optimization, making it well-suited for large-scale retail applications.

\subsubsection{Dynamic Programming}

Dynamic Programming (DP) is an optimization technique that decomposes a complex problem into a series of overlapping subproblems, solves each of them once, and reuses these solutions to efficiently construct an optimal solution to the original problem. By systematically building from smaller instances, DP avoids redundant computations and exploits the problem’s inherent optimal substructure and overlapping subproblem properties.
In the context of the linear binary knapsack problem, DP provides an exact and computationally efficient solution. Whereas brute-force enumeration explores all possible item combinations with exponential complexity, DP achieves optimality in pseudo-polynomial time by storing intermediate results in a structured manner. Specifically, a two-dimensional objective table $m[i,j]$ of size  $(n,W)$ is constructed, where each entry represents the maximum value attainable using the first 
$i$ items and a capacity limit $j$. If $w_i$ and $v_i$ denote the weight (space requirement) and value (objective contribution) of item $i$, respectively, the recurrence relation is defined as:

\begin{equation}
\resizebox{\columnwidth}{!}{$
m[i,j]=
\begin{cases}
m[i-1,j], & \text{if } w_i>j,\\
\max\big(m[i-1,j],\, m[i-1,j-w_i]+v_i\big), & \text{if } w_i\le j
\end{cases}
$}
\end{equation}
\vspace{-0.5em}

The computational complexity of the DP algorithm for the linear knapsack formulation is 
$\mathcal{O} (n c)$ where $n$ denotes the number of potential items in the POG and $c$
represents the total available capacity. In practice, the capacity $c$ can be scaled by adjusting the measurement units of shelf space (e.g., using half-inch or quarter-inch increments instead of finer units), thereby reducing the effective problem size without sacrificing solution quality. This makes DP both exact and scalable, providing a robust foundation for planogram-level optimization in large retail environments.

\subsection{Bay optimization}

This section develops and formulates the space optimization problem at the departmental level, where the objective is to determine the optimal allocation of available bays among different POGs within a given product category. For each category, the optimization solution specifies the amount of space—measured in bays—to be assigned to each POG, ensuring efficient utilization of the total available space.
The model takes as input several key parameters for each POG $i$, the minimum ($min_i$), the maximum ($max_i$) allowable bay allocations, the multiplicity constraint ($m_i$), which defines the smallest divisible unit of space (when fractional allocations are permissible), and the elasticity function $f_i$, represented by its coefficients, which describes the relationship between allocated space and expected sales (or margin).
Additionally, the total number of bays $s$ available for a category serves as a global constraint on space allocation. The elasticity function $f_i$ captures how sales or margin performance varies with changes in the allocated space for each POG. While the function need not always be represented explicitly by its coefficients, for illustrative purposes, logarithmic forms are employed in this section to model diminishing returns with respect to space expansion. The estimation of elasticity curves is performed by analyzing empirical dependencies between space and sales across store clusters, allowing the model to generalize performance patterns across similar retail environments.
However, accurately estimating these elasticity functions presents several challenges, particularly due to the impact of clustering on model stability and accuracy \cite{kok2007demand}. The composition and granularity of store clusters significantly influence the shape of elasticity curves and, consequently, the optimization results. Addressing these challenges motivated the development of Model Version 1 (v1), which incorporates improved clustering strategies and more robust estimation procedures for elasticity modeling.

\subsubsection{MIP formulation}
In this section we present a MIP for the problem. The integer variables for the formulation are the $x_i$. Together with the multiples $m_i$ they indicate the bay allocation to POG.

\begin{itemize}
\item[] $x_i$ indicates number of $m_i$ bays allocated the $i^\text{th}$ POG. For example, if $m_i = \frac{1}{2}$ then $x_i$ provides the number of half bays allocated to POG $i$. \end{itemize}

\begin{itemize}
\item[] $y_{i} = m_i x_i$ indicates number of bays allocated to the $i^\text{th}$ POG. If the multiple $m_i$ is not 1, then $y_i$ can assume fractional values.
\end{itemize}

The objective is to maximize the total number of sales (or margin) based on the bay allocation for each POG in the category:
$$\max \sum_{i=0}^{n} f_i(y_i) $$
Variable $y_i$ accounts for total bay allocation for all $i =0, \dots, n$:
\begin{equation}
\sum_{i=0}^{n} m_ix_i = y_i
\end{equation}
Each category has minimum and maximum bay allocation requirements:
\begin{equation}
min_i \leq y_i \leq max_i, \text{ for all } i =0, \dots, n
\end{equation}
We require that at most $s$ bays are allocated:
\begin{equation}
\sum_{i=0}^{n} m_ix_i \leq s
\end{equation}
The variables $x_i$ are non-negative integers:
\begin{equation}
x_i \in \mathbb{Z}, \text{ for all } i =0, \dots, n
\end{equation}
\begin{equation}
x_i \geq 0, \text{ for all } i =0, \dots, n
\end{equation}

\subsubsection{Solution}
In order to solve the above MIP formulation we linearize the functions $f_i$ to transform the formulation into a linear MIP.

Notice that the problem is separable because the decision variables $y_i$ appear separately in each function $f_i$. We use integer breakpoints between $min_i$ and $max_i$ to derive piecewise-linear approximation $f^a_i$ for $f_i$. For example, if $min_i = 2$ and $max_i = 4$ then the approximation curve $f^a_i(y_i) = f_i(2)\lambda_{i0} + f_i(2)\lambda_{i1}+f_i(2)\lambda_{i2}$, where

\begin{equation}
y_i = 2\lambda_{i0} + 3\lambda_{i1}+4\lambda_{i2}
\end{equation}

\begin{equation}
\lambda_{i0} + \lambda_{i1}+\lambda_{i2} = 1
\end{equation}

\begin{equation}
\lambda_{ij} \geq 0, \text{ for all } j = 0,1,2
\end{equation}

We assume that the objective function is concave. We conclude that the adjacency condition holds.

The transformed formulation in standard form is the following:
$$\max \sum_{i=0}^{n}  \sum_{j=min_i}^{max_i} \lambda_{ij}f_i(j) $$

\begin{equation}
z_i = - m_ix_i + \sum_{j=min_i}^{max_i} j\lambda_{ij}, \text{ for all } i = 0, \dots, n
\end{equation}

\begin{equation}
z_{n+k} =  \sum_{j=min_i}^{max_i} \lambda_{ij},  k =0, \dots, n
\end{equation}

\begin{equation}
\sum_{i=0}^{n} m_ix_i \leq s
\end{equation}

\begin{equation}
z_i = 0, \text{ for all } i =0, \dots, n
\end{equation}

\begin{equation}
z_{n+k} = 1, \text{ for all } k=0, \dots, n
\end{equation}

\begin{equation}
x_i \in \mathbb{Z}, \text{ for all } i =0, \dots, n
\end{equation}

\begin{equation}
x_i \geq 0, \text{ for all } i =0, \dots, n
\end{equation}

Upon constructing the MIP formulation, we solve it using OR-Tools.

\section{Results}
To evaluate the overall effectiveness of the OPTIMUS framework, we analyzed its end-to-end impact on sales and margin across ten recent optimization runs spanning multiple departments and store clusters. For each run, we compared pre- and post-optimization outcomes using historical data, capturing the combined effect of assortment selection and space allocation without isolating individual optimization stages.

Across the ten runs, OPTIMUS achieved an average sales lift of $11.83\% \pm 4.75\%$ and an average gross margin lift of $9.50\% \pm 3.27\%$. These results indicate consistent and substantial improvements in financial performance, with variability reflecting differences in department composition, elasticity characteristics, and baseline space utilization.

Overall, the findings demonstrate that OPTIMUS delivers robust, reproducible gains at the department level while remaining scalable across heterogeneous store environments. By jointly optimizing assortment and space decisions, the framework provides a practical and explainable decision-support tool for enterprise-scale retail space management.

\begin{table}[h!]
\centering
\caption{Overall performance lift across 10 optimization runs.}
\label{tab:overall_lift_pm}
\begin{tabular}{lc}
\toprule
\textbf{Metric} & \textbf{Lift (\%)} \\
\midrule
Sales  & $11.83 \pm 4.75$ \\
Margin & $9.50 \pm 3.27$ \\
\bottomrule
\end{tabular}
\end{table}

\section{Future Directions}
\label{sec:future}

While OPTIMUS currently uses a user-specified multi-objective weight vector
$\lambda = (\lambda_\text{sales}, \lambda_\text{margin}, \lambda_\text{units}, \lambda_\text{similarity})$,
a promising extension is to \emph{learn} these weights from historical outcomes rather than setting them manually.
Concretely, one can parameterize $\lambda$ as $\lambda_\theta \in \mathbb{R}^4_{+}$ (e.g., using a softplus transform)
and optimize $\theta$ to minimize a downstream loss defined on realized business outcomes (e.g., profit or margin lift,
with optional stability penalties for assortment churn).

A key challenge is that the planogram optimizer produces discrete decisions (a binary knapsack solution).
This can be addressed by using a differentiable surrogate during training, such as a continuous relaxation of the knapsack
decision variables, or perturbation-based gradient estimators, and then deploying the original exact dynamic programming
solver at inference time with the learned $\lambda$. This decision-focused learning direction would reduce manual tuning and
better align planogram decisions with end-to-end financial objectives, while preserving interpretability through the learned,
human-readable weights.

Additional future work includes incorporating substitution and cannibalization effects (moving beyond the linear knapsack
assumption), extending elasticity estimation with uncertainty-aware or robust optimization, and exploring multi-period
settings where planogram changes incur explicit operational switching costs.

\section{Conclusions}

We introduced OPTIMUS, an enterprise-scale framework for department-level retail space optimization that couples exact planogram assortment selection with departmental bay allocation. At the planogram level, OPTIMUS casts assortment construction as a linear binary knapsack problem and solves it via Dynamic Programming to produce deterministic, reproducible assortments under capacity constraints. These outputs feed a bay-level mixed-integer optimization model that allocates space across POGs while enforcing operational constraints. Across multiple departments and store clusters, the end-to-end results demonstrate that jointly optimizing \emph{what} to carry and \emph{how much space} to allocate can consistently improve financial performance while remaining scalable and interpretable for real-world merchandising workflows.

\bibliographystyle{icml2025}
\bibliography{ref}

\end{document}